\documentclass[11pt]{article}
\usepackage{mathtools}
\usepackage{amsfonts}
\usepackage{amsmath}
\usepackage{amsthm}
\usepackage{xfrac}
\usepackage{graphicx}
\usepackage{todonotes}
\usepackage{comment}
\usepackage{theoremref}
\usepackage{authblk}
\usepackage{citesort}

\newtheorem{lemma}{Lemma}
\newtheorem{theorem}{Theorem}

\begin{document}
\title{A Robust Preconditioner for High-Contrast Problems}
\author[1]{Yuliya Gorb\thanks{gorb@math.uh.edu, corresponding author}}
\author[1]{Daria Kurzanova\thanks{dariak@math.uh.edu}}
\author[1]{Yuri Kuznetsov\thanks{yuri@math.uh.edu}}
\affil[1]{Department of Mathematics,  University of Houston, Houston, TX 77204}
\date{}

\maketitle

\begin{abstract}
\noindent This paper concerns robust numerical treatment of an elliptic PDE with high contrast coefficients. 
A finite-element discretization of such an equation yields a linear system whose conditioning worsens as the variations in the values of PDE coefficients becomes large. This paper introduces a procedure by which the discrete system obtained from a linear finite element discretization of the given continuum problem is converted into an equivalent linear system of the saddle point type.  Then a robust preconditioner for the Lanczos method of minimized iterations for solving the derived saddle point problem is proposed. Robustness with respect to the contrast parameter and the mesh size is justified. Numerical examples support theoretical results and demonstrate independence of the number of iterations on the contrast, the mesh size and also on the different contrasts on the inclusions. 
 \end{abstract}


\section{Introduction}

\begin{figure}[ht]
\centering
\includegraphics[scale=.75]{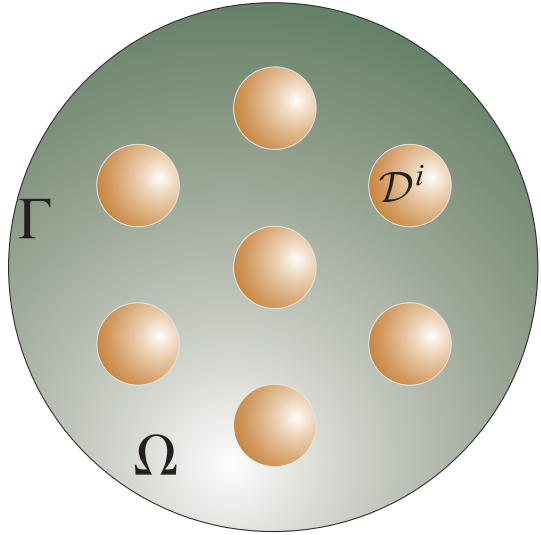} \caption{The domain $\Omega$ with highly conducting inclusions $\mathcal{D}^i$, $i\in \{1,\ldots,m\}$} \label{F:main}
\end{figure}

In this paper, we consider an iterative solution of the linear system arising from the discretization of the diffusion problem
\begin{equation} \label{E:problem-intro}
-\nabla \cdot \left[  \sigma(x) \nabla u \right]  = f, \quad x\in  \Omega
\end{equation}
with appropriate boundary conditions on $\Gamma = \partial \Omega$. We assume that $\Omega$ is a bounded domain $\Omega \subset \mathbb{R}^d$, $d\in \{2,3\}$, that contains $m\geq 1$ polygonal or polyhedral subdomains $\mathcal{D}^i$, see Fig. \ref{F:main}. 
Also assume that the distance between the neighboring $\mathcal{D}^i$ and $\mathcal{D}^j$ is at least of order of the sizes of these subdomains, that is, bounded below by a multiple of their diameters.
The main focus of this work is on the case when the coefficient function $ \sigma(x) \in L^{\infty} (\Omega)$ varies largely within the domain $\Omega$, that is,
\[
\kappa = \frac{\sup_{x\in \Omega}  \sigma(x)}{\inf_{x\in \Omega}  \sigma(x)} \gg 1.
\]
In this work, we assume that the domain $\Omega$ contains disjoint polygonal or polyhedral subdomains $\mathcal{D}^i$, $i \in \{1,\ldots,m\}$, where $\sigma$ takes ``large'' values, e.g. of order $O(\kappa)$, but remains of $O(1)$ in the domain outside of $\mathcal{D} := \cup_{i=1}^m \mathcal{D}^i$. 

The P1-FEM discretization of this problem results in a linear system
\begin{equation} \label{E:linsys}
\boldsymbol{\mathcal{K}} \overline{u}  = \overline{F},
\end{equation}
with a large and sparse matrix $\boldsymbol{\mathcal{K}}$. 
A major issue in numerical treatments of \eqref{E:problem-intro}, with the coefficient $\sigma$ discussed above, is that the high contrast leads to an ill-conditioned matrix $\boldsymbol{\mathcal{K}}$ in \eqref{E:linsys}. Indeed, if $h$ is the discretization scale, then the condition number of the resulting stiffness matrix $\boldsymbol{\mathcal{K}}$ grows proportionally to $h^{-2}$ with coefficient of proportionality depending on $\kappa$. Because of that, the high contrast problems have been a subject of an active research recently, see e.g. \cite{ah02,agks08}. 

There is one more feature of the system \eqref{E:linsys} that we investigate in this paper. Recall that 
if $\boldsymbol{\mathcal{K}}$ is symmetric and positive definite, then \eqref{E:linsys} is typically solved with the Conjugate Gradient (CG) method, if $\boldsymbol{\mathcal{K}}$ is nonsymmetric then the most common solver for \eqref{E:linsys} is  GMRES. In this paper, we introduce an additional variable that allows us to replace \eqref{E:linsys} with an equivalent formulation in terms of a linear system 
 \begin{equation} \label{E:LS}
\boldsymbol{\mathcal{A}}\boldsymbol{x}=\boldsymbol{\mathcal{F}}
 \end{equation}
 with a {\it saddle point matrix}  $\boldsymbol{\mathcal{A}}$ written in the block form:
 \begin{equation} \label{E:Ablock}
\boldsymbol{\mathcal{A}}=
\begin{bmatrix}
\bold{A} & \bold{B}^T \\
\bold{B} & -\bold{\Sigma}
\end{bmatrix},
 \end{equation}
where $\bold{A} \in \mathbb{R}^{n \times n}$ is symmetric and positive definite, $\bold{B}  \in \mathbb{R}^{k \times n} $ is
rank deficient, and $\bold{\Sigma} \in \mathbb{R}^{k \times k} $ is symmetric and positive semidefinite, so that the corresponding linear system \eqref{E:LS} is singular but consistent. 
Unfortunately, Krylov space iterative methods tend to converge very slowly when applied to systems with saddle point matrices and preconditioners are needed to achieve faster convergence. 
The CG method that was mainly developed for the iterative solution of linear systems with symmetric and definite matrices is not, in general, robust for systems with indefinite matrices, \cite{wsy04}. The {\it Lanczos algorithm} of minimized iterations does not have such a restriction and has been utilized in this paper. 
Below, we introduce a construction of a robust preconditioner for solving \eqref{E:LS} by the Lanczos iterative scheme \cite{km74}, that is, whose convergence rate is independent of the contrast parameter $\kappa \gg 1$ and the discretization size $h>0$.  

Also, the special case of \eqref{E:linsys} with \eqref{E:Ablock} considered in the Appendix of this paper is when $\bold{\Sigma} \equiv \bold{0}$. The problem of this type has received considerable attention over the years. But the most studied case is when $\boldsymbol{\mathcal{A}}$ is {\it nonsingular}, in which case $\bold{B}$ must be of full rank, see e.g. \cite{kgw00,lv98} and references therein. The main focus of this paper is on singular $\boldsymbol{\mathcal{A}}$ with the rank deficient block $\bold{B}$. Below, we propose a block-diagonal preconditioner for the Lanczos method employed to solve the problem  \eqref{E:LS}, and this preconditioner is also singular. We also rigorously justify its robustness with respect to $h$ and $\kappa$. Our numerical experiments on simple test cases support our theoretical findings.

Finally, we point out that a robust numerical treatment of the described problem
is crucial in developing the mutiscale strategies for models of composite materials with 
highly conducting particles. The latter find their application in particulate flows, subsurface flows in natural porous formations,
electrical conduction in composite materials, and medical and geophysical imaging.

\vspace{5pt}

The rest of this paper is organized as follows. In Chapter \ref{S:form-results} the mathematical problem formulation is presented and main results are stated. Chapter \ref{S:proof} discusses proofs of main results, and numerical results of the proposed procedure are given in Chapter \ref{S:numercs}. Conclusions are presented in Chapter \ref{S:concl}. Proof an auxiliary fact is given in Appendix. 

\vspace{8pt}

\noindent {\bf Acknowledgements.} 
First two authors were supported by the NSF grant DMS-$1350248$.

\section{Problem Formulation and Main Results} \label{S:form-results}
Consider an open, a bounded domain $\Omega \subset \mathbb{R}^d$, $d\in \{2,3\}$ with piece-wise smooth boundary $\Gamma$, that contains $m\geq 1$ subdomains $\mathcal{D}^i$, which are located at distances comparable to their sizes from one another, see Fig. \ref{F:main}. For simplicity, we assume that $\Omega$ and $\mathcal{D}^i$ are polygons if $d=2$ or polyhedra if $d=3$. The union of $\mathcal{D}^i$ is denoted by $\mathcal{D}$.
In the domain $\Omega$ we consider the following elliptic problem
\begin{equation} \label{E:pde-form}
\left\{
\begin{array}{r l l}
-\nabla \cdot \left[  \sigma(x) \nabla u \right] & = f, & x\in  \Omega  \\[2pt]
u & = 0, & x\in  \Gamma
\end{array}
\right.
\end{equation}
with the coefficient $\sigma$ that largely varies inside the domain $\Omega$. 
For simplicity of the presentation, we focus on the case when $\sigma$ is a piecewise constant function given by
\begin{equation} \label{E:sigma}
\sigma(x)=
\begin{cases}
1, & x\in  \Omega\setminus \overline{\mathcal{D}}\\[5pt]
\displaystyle 1+\frac{1}{\varepsilon_i}, &  x\in \mathcal{D}_i, ~i\in\{1,\ldots,m\}
\end{cases}
\end{equation}
with $\displaystyle \varepsilon:=\max_i \varepsilon_i \ll 1$. We also assume that the source term in \eqref{E:pde-form} is $f \in L^2(\Omega)$. 

When performing a P1-FEM discretization of \eqref{E:pde-form} with \eqref{E:sigma}, we choose a FEM space $V_h \subset H_0^1 (\Omega)$ to be the space of linear finite-element functions defined on a conforming
quasi-uniform triangulation $\Omega_h$ of $\Omega$ of the size $h\ll 1$. For simplicity, we assume that $\partial \Omega_h = \Gamma$. With that, the classical FEM discretization results in the system of the type \eqref{E:linsys}. We proceed differently and derive another discretized system of the saddle point type as shown below.

\subsection{Derivation of a Singular Saddle Point Problem}

If $\left.  \mathcal{D}^i_h = \Omega_h\right|_{\mathcal{D}^i}$  then we denote $\left. V_h^i : = V_h\right|_{\mathcal{D}^i_h}$ and $ \mathcal{D}_h := \cup_{i=1}^m   \mathcal{D}^i_h$. With that, we write the FEM formulation of  \eqref{E:pde-form}-\eqref{E:sigma} as
\[ \mbox{Find}\quad u_h \in V_h \quad \mbox{and} \quad \lambda_h = (\lambda^1_h,\ldots,\lambda^m_h) \quad \mbox{with}\quad  \lambda^i_h \in V^i_h \quad \mbox{such that}
\] 
\begin{equation} \label{E:FEM-form}
\int_{\Omega_h} \nabla u_h \cdot  \nabla v_h ~dx + \int_{\mathcal{D}_h} \nabla \lambda_h \cdot  \nabla v_h ~dx = \int_{\Omega_h} f v_h ~dx, \quad \forall v_h\in V_h,
\end{equation}
provided 
\begin{equation} \label{E:FEM-cond}
u_h =\varepsilon_i \lambda^i_h+c_i  \quad \mbox{in} \quad\mathcal{D}^i_h, \quad i\in \{1,\ldots,m\},
\end{equation}
where $c_i$ is an arbitrary constant. 
First, we turn out attention to the FEM discretization of \eqref{E:FEM-form} that yields a system of linear equations
\begin{equation} \label{E:lin-sys}
\bold{A} \overline{u} + \bold{B}^T \overline{\lambda}  =\overline{ \mathrm{F}},
\end{equation}
and then discuss implications of \eqref{E:FEM-cond}.

To provide the comprehensive description of all elements of the system \eqref{E:lin-sys}, we introduce the following notations for the number of degrees of freedom in different parts of  $\Omega_h$. Let $N$ be the total number of nodes in $\Omega_h$, and 
$n$ be the number of nodes in $\overline{\mathcal{D}}_h$ so that 
\[
n=\sum_{i=1}^m n_i,
\]
where $n_i$ denotes the number of degrees of freedom in  $\overline{\mathcal{D}}^i_h$, and, finally, $n_{0}$ is the number of nodes in $\Omega_h \setminus \overline{\mathcal{D}_h}$, so that we have
\[
N=n_{0} + n = n_{0}+\sum_{i=1}^m n_i .
\]
Then in \eqref{E:lin-sys}, the vector $\overline{u} \in \mathbb{R}^{N}$ has entries $u_i=u_h(x_i) $ with $x_i \in \overline{\Omega}_h$. We count the entries of $\overline{u}$ in such a way that its first $n$ elements correspond to the nodes of $\overline{\mathcal{D}}_h$, and the remaining $n_0$ entries correspond to the nodes of $\overline{\Omega}_h \setminus \overline{ \mathcal{D}}_h$. 
Similarly, the vector $\overline{\lambda}\in \mathbb{R}^{n}$ has entries $\lambda_i=\lambda_h(x_i) $ where $x_i \in \overline{ \mathcal{D}}_h$. 

The symmetric positive definite matrix $\bold{A} \in\mathbb{R}^{N\times N}$ of \eqref{E:lin-sys}  is the stiffness matrix that arises from the discretization of the Laplace operator with the homogeneous Dirichlet boundary conditions on $\Gamma$. 
Entries of $\bold{A}$ are defined by
\begin{equation} \label{E:A-def}
(\bold{A} \overline{u} , \overline{v} ) = \int_{\Omega_h} \nabla u_h \cdot \nabla v_h ~dx, \quad\mbox{where}\quad \overline{u} , \overline{v} \in \mathbb{R}^N, \quad
u_h , v_h \in V_h,
\end{equation}
where $(\cdot,\cdot)$ is the standard dot-product of vectors.
This matrix can also be partitioned into
\begin{equation} \label{E:matr-A}
\bold{A} = 
\begin{bmatrix}
\mathrm{A}_{\mathcal{D}\mathcal{D}} & \mathrm{A}_{\mathcal{D}0} \\
\mathrm{A}_{0\mathcal{D}} & \mathrm{A}_{00} 
\end{bmatrix},
\end{equation}
where the block $\mathrm{A}_{\mathcal{D}\mathcal{D}} \in  \mathbb{R}^{n\times n}$ is the stiffness matrix corresponding to the highly conducting  inclusions $\overline{\mathcal{D}}^i_h$, $i\in \{1,\ldots,m\}$, the  block $\mathrm{A}_{00}\in  \mathbb{R}^{n_0\times n_0}$ corresponds to the region outside of $\overline{\mathcal{D}}_h$, and the entries of $\mathrm{A}_{\mathcal{D}0} \in  \mathbb{R}^{n\times n_0}$ and $\mathrm{A}_{0\mathcal{D}}=\mathrm{A}_{\mathcal{D}0}^T$ are assembled from contributions both from finite elements in $\overline{\mathcal{D}}_h$ and $\overline{\Omega}_h \setminus \overline{\mathcal{D}}_h$.

The matrix $\bold{B} \in \mathbb{R}^{n\times N}$ of \eqref{E:lin-sys} is also written in the block form  as
\begin{equation} \label{E:matr-B}
\bold{B} = \begin{bmatrix} \boldsymbol{\mathcal{B}}_{\mathcal{D}} & \bold{0}\end{bmatrix}
\end{equation}
with zero-matrix $\bold{0} \in \mathbb{R}^{n\times n_0}$ and $\boldsymbol{\mathcal{B}}_{\mathcal{D}} \in \mathbb{R}^{n\times n}$ that corresponds to the highly conducting inclusions. The matrix $\boldsymbol{\mathcal{B}}_{\mathcal{D}}$ is the stiffness matrix in the discretization of the Laplace operator in the domain $\overline{\mathcal{D}}_h$ with the Neumann boundary conditions on $\partial \mathcal{D}_h$. 
In its turn, $\boldsymbol{\mathcal{B}}_{\mathcal{D}}$ is written in the block form as
\[
\boldsymbol{\mathcal{B}}_{\mathcal{D}} = 
\begin{bmatrix}
\mathrm{B}_{1} & \ldots & 0\\
\vdots  &\ddots & \vdots \\
0 & \ldots & \mathrm{B}_{m} 
\end{bmatrix}=\text{diag}~ (\mathrm{B}_{1}, \ldots, \mathrm{B}_{m}),
\]
with matrices $\mathrm{B}_i$$\in \mathbb{R}^{n_i \times n_i}$, whose entries are similarly defined by
\begin{equation} \label{E:B-def}
(\mathrm{B}_i \overline{u} , \overline{v} ) = \int_{\mathcal{D}^i_h} \nabla u_h \cdot \nabla v_h ~dx, \quad\mbox{where}\quad \overline{u} , \overline{v} \in \mathbb{R}^{n_i}, \quad
u_h , v_h \in V^i_h.
\end{equation}
We remark that each $\mathrm{B}_i$ is positive semidefinite with 
\begin{equation} \label{E:ker-Bi}
\ker \mathrm{B}_i = \mbox{span}
\left\{\begin{bmatrix}
1 \\
\vdots\\
1
\end{bmatrix}\right \}.
\end{equation}

Finally, the vector $\overline{ \mathrm{F}}  \in \mathbb{R}^{N}$ of  \eqref{E:lin-sys}  is defined in a similar way by
\[
(\overline{ \mathrm{F}},\overline{v}) = \int_{\Omega_h} f v_h ~dx, \quad\mbox{where}\quad \overline{v} \in \mathbb{R}^N, \quad v_h \in V_h.
\]

To complete the derivation of the linear system corresponding to \eqref{E:FEM-form}-\eqref{E:FEM-cond}, we rewrite \eqref{E:FEM-cond} in the weak form that is as follows:
\begin{equation} \label{E:FEM-cond-new}
 \int_{\mathcal{D}^i_h} \nabla u_h \cdot  \nabla v^i_h ~dx -  \varepsilon_i \int_{\mathcal{D}^i_h} \nabla \lambda^i_h 
 \cdot  \nabla v^i_h ~dx=0 , \quad i\in\{1,\ldots,m\} \quad \mbox{for all } v^i_h \in V^i_h , 
\end{equation}
and add the discrete analog of  \eqref{E:FEM-cond-new} to the system \eqref{E:lin-sys}. For that, denote
\[
\bold{\Sigma}_\varepsilon =  
\begin{bmatrix}
\varepsilon_1\mathrm{B}_{1} & \ldots & 0\\
\vdots  &\ddots & \vdots \\
0 & \ldots & \varepsilon_m\mathrm{B}_{m} 
\end{bmatrix}=\text{diag}~ (\varepsilon_1\mathrm{B}_{1}, \ldots, \varepsilon_m\mathrm{B}_{m}),
\]
then \eqref{E:FEM-cond-new} implies
\begin{equation} \label{E:lambda-u-relation}
\bold{B} \overline{u} - \bold{\Sigma}_\varepsilon  \overline{\lambda}  =  \overline{0}.
\end{equation}
This together with \eqref{E:lin-sys} yields
\begin{equation} \label{E:lin-sys-full}
\left\{
\begin{array}{r l l}
\bold{A} \overline{u} + \bold{B}^T \overline{\lambda} & =\overline{ \mathrm{F}}, \\[2pt]
\bold{B} \overline{u} - \bold{\Sigma}_\varepsilon \overline{\lambda}  & =\overline{0},
\end{array}
\right.\quad \overline{u} \in \mathbb{R}^N, \quad  \overline{\lambda} \in \mathbb{R}^n, ~
\overline{\lambda} \in \textrm{Im }  \boldsymbol{\mathcal{B}}_{\mathcal{D}},
\end{equation}
or
\begin{equation} \label{E:lin-sys-full-matrix-1}
 \boldsymbol{\mathcal{A}}_\varepsilon \bold{x}_\varepsilon = \overline{ \mathcal{F}} ,
\end{equation}
where
\begin{equation} \label{E:lin-sys-full-matrix-2}
\boldsymbol{\mathcal{A}}_\varepsilon = \begin{bmatrix} \bold{A} & \bold{B}^T \\ \bold{B} & - \bold{\Sigma}_\varepsilon  \end{bmatrix}
= \begin{bmatrix} \mathrm{A}_{\mathcal{D}\mathcal{D}} & \mathrm{A}_{\mathcal{D}0} & \boldsymbol{\mathcal{B}}_{\mathcal{D}} \\ \mathrm{A}_{0\mathcal{D}} & \mathrm{A}_{00} &  \bold{0}^T \\ \boldsymbol{\mathcal{B}}_{\mathcal{D}}  & \bold{0} & -\bold{\Sigma}_\varepsilon \end{bmatrix}, \quad 
\bold{x}_\varepsilon = 
\begin{bmatrix}
\overline{u} \\
\overline{\lambda} 
\end{bmatrix},
\quad  \overline{ \mathcal{F}}  = \begin{bmatrix} \overline{ \mathrm{F}} \\ \overline{0} \end{bmatrix}.
\end{equation}
This saddle point formulation \eqref{E:lin-sys-full-matrix-1}-\eqref{E:lin-sys-full-matrix-2} for the PDE \eqref{E:pde-form}-\eqref{E:sigma} was first proposed in \cite{kuz09}.
Clearly, there exists a unique solution $\overline{u} \in \mathbb{R}^N$ and $\overline{\lambda} \in \mathbb{R}^n$, $\overline{\lambda} \in \textrm{Im }  \boldsymbol{\mathcal{B}}_{\mathcal{D}}$ of \eqref{E:lin-sys-full-matrix-1}-\eqref{E:lin-sys-full-matrix-2}. 

It is important to point out that the main feature of the problem \eqref{E:lin-sys-full} is in rank deficiency of the matrix $\bold{B}$. This would lead to the introduction in the next Section \ref{S:main-results} of a {\it singular} block-diagonal preconditioner for the Lanczos method employed to solve the problem  \eqref{E:lin-sys-full-matrix-2}. 
Independence of the convergence of the employed Lanczos method on the discretization size $h>0$ follows from the spectral properties of the constructed preconditioner that are independent of $h$ due to the norm preserving extension theorem of \cite{widlund87}. Independence on contrast parameters $\varepsilon_i$ follows from the closeness of spectral properties of the matrices of the original system 
\eqref{E:lin-sys-full-matrix-2} and the limiting one \eqref{E:lin-sys-infinite-matrix}, also demonstrated in Appendix. Our numerical experiments below also show independence of the iterative procedure on the number of different contrasts $\varepsilon_i$, $i \in \{1,\ldots,m\}$, in the  inclusions $\mathcal{D}^i$.

\subsection{Preconditioned System and Its Implementation}

\subsubsection{Lanczos Method}

In principal, we could have used the CG method that was mainly developed for the iterative solution of linear systems with symmetric and {\it definite} matrices, and apply it to the square of the matrix of the preconditioned system. However, the {\it Lanczos method} of minimized iterations is not restricted to the definite matrices, and, since it has the same arithmetic cost as CG, is employed in this paper. A symmetric and positive semidefinite block-diagonal preconditioner of the form
\begin{equation}   \label{E:block-prec}
\boldsymbol{\mathcal{P}} = \begin{bmatrix}
\mathcal{P}_\mathrm{A} & 0  \\
0 & \mathcal{P}_{\mathrm{B}}
\end{bmatrix},
\end{equation}
for this method is also proposed in this section, where the role of the blocks $\mathcal{P}_\mathrm{A} $ and $\mathcal{P}_\mathrm{B}$ will be explained below. But, first, for the completeness of presentation, we describe the Lanczos algorithm.

For a symmetric and positive semidefinite matrix $\boldsymbol{\mathcal{H}}$ that later will be defined as the Moore-Penrose 
pseudo inverse\footnote{$\mathrm{M}^{\dagger}$ is the Moore-Penrose pseudo inverse of $\mathrm{M}$ if and only if it satisfies the following Moore-Penrose equations, see e.g.  \cite{axel}:
\[
\mathrm{(i)}~\mathrm{M}^{\dagger}\mathrm{M}\mathrm{M}^{\dagger}=\mathrm{M}^{\dagger}, \quad 
\mathrm{(ii)}~\mathrm{M}\mathrm{M}^{\dagger}\mathrm{M} =\mathrm{M},   \quad 
\mathrm{(iii)} ~\mathrm{M}\mathrm{M}^{\dagger}~ \mbox{and}~\mathrm{M}^{\dagger}\mathrm{M}~ \mbox{are symmetric.}
\]} $ \boldsymbol{\mathcal{P}}^{\dagger} $, 
introduce a new scalar product
\[
(\overline{x},\overline{y})_{\mathcal{H}} := (\boldsymbol{\mathcal{H}}\overline{x},\overline{y}), \quad  \mbox{for all} \quad \overline{x},\overline{y} \in \mathbb{R}^{N+n},~ \overline{x},~\overline{y}\perp \ker{\boldsymbol{\mathcal{H}}},
\]
and consider the preconditioned Lanczos iterations, see \cite{km74}, $\overline{z}^k = \begin{bmatrix} \overline{u}^k\\ \overline{\lambda}^k \end{bmatrix} \in \mathbb{R}^{N + n}$, $k\geq 1$:
\[
\overline{z}^k = \overline{z}^{k-1} - \beta_k \overline{y}_k,  
\]
where
\[
\beta_k = \frac{( \boldsymbol{\mathcal{A}}_{\varepsilon}\overline{z}^{k-1}-\overline{\mathcal{F}}, \boldsymbol{\mathcal{A}}_{\varepsilon} \overline{y}_{k})_{\mathcal{H}}}{(\boldsymbol{\mathcal{A}}_{\varepsilon} \overline{y}_{k}, \boldsymbol{\mathcal{A}}_{\varepsilon} \overline{y}_{k})_{ \mathcal{H}}} .
\]
and
\[
y_k = \begin{cases} 
\boldsymbol{\mathcal{H}}(\boldsymbol{\mathcal{A}}_{\varepsilon}\overline{z}^0 - \overline{\mathcal{F}}), &k=1 \\
\boldsymbol{\mathcal{H}}\boldsymbol{\mathcal{A}}_{\varepsilon}\overline{y}_1 - \alpha_2 \overline{y}_1, & k=2\\
\boldsymbol{\mathcal{H}}\boldsymbol{\mathcal{A}}_{\varepsilon}\overline{y}_{k-1} - \alpha_k \overline{y}_{k-1} - \gamma_k \overline{y}_{k-2}, & k>2,
\end{cases}
\]
with
\[
 \alpha_k = \frac{(\boldsymbol{\mathcal{A}}_{\varepsilon} \boldsymbol{\mathcal{H}} \boldsymbol{\mathcal{A}}_{\varepsilon} \overline{y}_{k-1}, \boldsymbol{\mathcal{A}}_{\varepsilon} \overline{y}_{k-1})_{ \mathcal{H}}}{(\boldsymbol{\mathcal{A}}_{\varepsilon} \overline{y}_{k-1}, \boldsymbol{\mathcal{A}}_{\varepsilon} \overline{y}_{k-1})_{\mathcal{H}}}, \qquad  
 \gamma_k = \frac{(\boldsymbol{\mathcal{A}}_{\varepsilon} \boldsymbol{\mathcal{H}}\boldsymbol{\mathcal{A}}_{\varepsilon} \overline{y}_{k-1}, \boldsymbol{\mathcal{A}}_{\varepsilon} \overline{y}_{k-1})_{ \mathcal{H}}}{(\boldsymbol{\mathcal{A}}_{\varepsilon} \overline{y}_{k-2}, \boldsymbol{\mathcal{A}}_{\varepsilon} \overline{y}_{k-2})_{ \mathcal{H}}}.
\]

\subsubsection{Proposed Preconditioner}

It was previously observed, see e.g. \cite{kuz09,kuz95}, that the following matrix
\begin{equation}   \label{E:theor-precond}
\bold{P}= \begin{bmatrix}
\bold{A} & \bold{0}  \\
\bold{0} & \bold{B}\bold{A}^{-1}\bold{B}^T
\end{bmatrix},
\end{equation} 
is the best choice for a block-diagonal preconditioner of $\boldsymbol{\mathcal{A}}_\varepsilon$. 
This is because the eigenvalues of the generalized eigenvalue problem
\begin{equation}   \label{E:gen-eigenpr-eps}
\boldsymbol{\mathcal{A}}_\varepsilon \bold{x}  = \mu \bold{P} \bold{x}, \quad \overline{u} \in \mathbb{R}^N, \quad \overline{\lambda} \in \mathrm{Im}\,\boldsymbol{\mathcal{B}}_{\mathcal{D}} ,
\end{equation} 
belong to the union of $[c_1,c_2]\cup[c_3,c_4]$ with $c_1\leq c_2 <0$ and $0<c_3\leq c_4$, with numbers $c_i$ being independent of both $h$, and $\varepsilon_i$, see \cite{irt93,kuz95,kuz09}. For the reader's convenience, the proof of this statement is also shown in Appendix below (see Lemma 6). 

The preconditioner $\bold{P}$ of \eqref{E:theor-precond} is of limited practical use and is a subject of primarily theoretical interest.  To construct a preconditioner that one can actually use in practice, we will find a matrix $\boldsymbol{\mathcal{P}}$ such that there exist constants $\alpha$, $\beta$ independent on the mesh size $h$ and that 
\begin{equation}   \label{E:spect-eq}
\alpha ( \bold{P} \bold{x}  , \bold{x}  ) \leq (\boldsymbol{\mathcal{P}}  \bold{x}  , \bold{x}   ) \leq \beta ( \bold{P}\bold{x}  , \bold{x}  )\quad \mbox{for all } \bold{x}  \in \mathbb{R}^{N+n}.
\end{equation} 
This property \eqref{E:spect-eq} is hereafter referred to as  {\it spectral equivalence} of $\boldsymbol{\mathcal{P}} $ to $\bold{P}$ of \eqref{E:theor-precond}. 
Obviously, the matrix $\boldsymbol{\mathcal{P}}$ of the form
\eqref{E:block-prec} has to be such that the block $\mathcal{P}_\mathrm{A} $ is spectrally equivalent to $\bold{A}$, whereas  $\mathcal{P}_{\mathrm{B}}$ is spectrally equivalent to $\bold{B}\bold{A}^{-1}\bold{B}^T$, see also \cite{irt93,kuz95,kuz09}.

For the block $\mathcal{P}_\mathrm{A} $, one can use any existing symmetric and positive definite preconditioner devised for the discrete Laplace operator on quasi-uniform and regular meshes. Note that for a regular hierarchical mesh, the best preconditioner for $\bold{A}$ would be the BPX preconditioner, see \cite{bpx90}. However, to extend our results to the hierarchical meshes, one needs the corresponding norm preserving extension theorem as in \cite{widlund87}. 
Hence, this paper is not investigating the effect of the choice $\mathcal{P}_\mathrm{A} $, and our primary aim is to propose a preconditioner $\mathcal{P}_{\mathrm{B}}$ that could be effectively used in solving \eqref{E:lin-sys-full}. 

To that end, for our  Lanczos method of minimized iterations, we use the following block-diagonal preconditioner:
\begin{equation} \label{E:our-prec}
\boldsymbol{\mathcal{P}} =
\begin{bmatrix}
\mathcal{P}_{\mathrm{A}} & 0\\
0 & \boldsymbol{\mathcal{B}}_{\mathcal{D}}
\end{bmatrix},
\end{equation} 
and in the numerical experiments below, we will simply take $\mathcal{P}_{\mathrm{A}}=\bold{A}$. 
Finally, we define
\begin{equation} \label{E:our-H}
\boldsymbol{\mathcal{H}} = \boldsymbol{\mathcal{P}}^{\dagger} =
\begin{bmatrix}
\mathcal{P}_{\mathrm{A}}^{-1} & 0\\
0 & \left[\boldsymbol{\mathcal{B}}_{\mathcal{D}}\right]^{\dagger}
\end{bmatrix},
\end{equation}
and remark that even though the matrix $\boldsymbol{\mathcal{B}}_{\mathcal{D}}$ is singular, as evident from the Lanczos algorithm above, one actually never needs to use its pseudo inverse at all. Indeed, this is due to the block-diagonal structure \eqref{E:our-H} of $\boldsymbol{\mathcal{H}}$, and the block form \eqref{E:lin-sys-full-matrix-2} of the original matrix $\boldsymbol{\mathcal{A}}_{\varepsilon}$.

\subsection{Main Result: Block-Diagonal Preconditioner} \label{S:main-results}

The main theoretical result of this paper establishes a robust preconditioner for solving \eqref{E:lin-sys-infinite-matrix} or, equivalently \eqref{E:lin-sys-infinite}, and is given in the following theorem.

\begin{theorem} \label{T:main}
Let the triangulation $\Omega_h$ for \eqref{E:pde-form}-\eqref{E:sigma} be conforming and quasi-uniform. 
Then the matrix $\boldsymbol{\mathcal{B}}_{\mathcal{D}}$ is spectrally equivalent to the matrix $\bold{B}\bold{A}^{-1}\bold{B}^T$, that is, 
there exist constants $\mu_\star, \mu^\star >0$ independent of $h$ and such that
\begin{equation} \label{E:main-thm}
\mu_\star \leq  \frac{  \left( \boldsymbol{\mathcal{B}}_{\mathcal{D}}\overline{\psi},\overline{\psi} \right)  }{  \left(  \bold{B}\bold{A}^{-1}\bold{B}^T \overline{\psi},\overline{\psi} \right)   } \leq \mu^\star,  \quad \mbox{for all}\quad  0\neq \overline{\psi}\in \mathbb{R}^n, 
~ \overline{\psi} \in \mathrm{Im}\,\boldsymbol{\mathcal{B}}_{\mathcal{D}}.
\end{equation}
\end{theorem}

\noindent This theorem asserts that the nonzero eigenvalues of the generalized eigenproblem
\begin{equation} \label{E:eigprob-1}
\bold{B}\bold{A}^{-1}\bold{B}^T \overline{\psi} = \mu \boldsymbol{\mathcal{B}}_{\mathcal{D}} \overline{\psi}, \quad \overline{\psi}\in \mathbb{R}^n,~
\overline{\psi} \in \mathrm{Im}\,\boldsymbol{\mathcal{B}}_{\mathcal{D}},
\end{equation}
are bounded. Hence, its proof is based on the construction of the {\bf upper} and {\bf lower} bounds for $\mu$ in \eqref{E:eigprob-1}
and is comprised of the following facts many of which are proven in the next section.

\begin{lemma} \label{L:fact1}
The following equality of matrices holds
\begin{equation} \label{E:fact1}
\bold{B}\bold{A}^{-1}\bold{B}^T = \boldsymbol{\mathcal{B}}_{\mathcal{D}} \mathrm{S}_{00}^{-1}  \boldsymbol{\mathcal{B}}_{\mathcal{D}},
\end{equation}
where
\[
\mathrm{S}_{00} = \mathrm{A}_{\mathcal{DD}} - \mathrm{A}_{\mathcal{D}0} \mathrm{A}^{-1}_{00}\mathrm{A}_{0\mathcal{D}},
\]
is the Schur complement to the block $A_{00}$ of the matrix $\bold{A}$ of \eqref{E:lin-sys-infinite-matrix}.
\end{lemma} 
\noindent This fact is straightforward and comes from the block structure of matrices $\bold{A}$ of \eqref{E:matr-A} and $\bold{B}$ of \eqref{E:matr-B}. Indeed, 
using this, the generalized eigenproblem \eqref{E:eigprob-1} can be rewritten as 
\begin{equation} \label{E:eigprob-2}
\boldsymbol{\mathcal{B}}_{\mathcal{D}}  \mathrm{S}_{00}^{-1} \boldsymbol{\mathcal{B}}_{\mathcal{D}}  \, \overline{\psi} = \mu \boldsymbol{\mathcal{B}}_{\mathcal{D}} \overline{\psi}, \quad  \quad \overline{\psi} \in  \mathbb{R}^{n},
~ \overline{\psi} \in \mathrm{Im }\,\boldsymbol{\mathcal{B}}_{\mathcal{D}}.
\end{equation}
Introduce a matrix $\mathrm{B}^{\sfrac{1}{2}}_{\mathcal{D}}$ via $\boldsymbol{\mathcal{B}}_{\mathcal{D}} =\mathrm{B}^{\sfrac{1}{2}}_{\mathcal{D}}\mathrm{B}^{\sfrac{1}{2}}_{\mathcal{D}}$ and note that $\ker \boldsymbol{\mathcal{B}}_{\mathcal{D}}  = \ker \mathrm{B}^{\sfrac{1}{2}}_{\mathcal{D}}$.

\begin{lemma} \label{L:fact2} 
The generalized eigenvalue problem \eqref{E:eigprob-2} is equivalent to 
\begin{equation} \label{E:eigprob-3}
\boldsymbol{\mathcal{B}}_{\mathcal{D}} ^{\sfrac{1}{2}} \mathrm{S}_{00}^{-1} \boldsymbol{\mathcal{B}}_{\mathcal{D}} ^{\sfrac{1}{2}}  \,\overline{\varphi} = \mu\,\overline{\varphi},
\quad  \overline{\varphi} \in \mathbb{R}^{n},~\overline{\varphi}  \in \mathrm{Im }\,\boldsymbol{\mathcal{B}}_{\mathcal{D}} ,
\end{equation}
in the sense that they both have the same eigenvalues $\mu$'s, and the corresponding eigenvectors are related via $\overline{\varphi} = \boldsymbol{\mathcal{B}}_{\mathcal{D}}^{\sfrac{1}{2}}\overline{\psi} \in \mathrm{Im}\,\boldsymbol{\mathcal{B}}_{\mathcal{D}}$.
\end{lemma} 

\begin{lemma} \label{L:fact3}
The  generalized  eigenvalue problem \eqref{E:eigprob-3} is equivalent to 
\begin{equation} \label{E:eigprob-4}
\boldsymbol{\mathcal{B}}_{\mathcal{D}} \,  \overline{u}_{\mathcal{D}}  = \mu \mathrm{S}_{00} \, \overline{u}_{\mathcal{D}},
\quad \overline{u}_{\mathcal{D}} \in \mathbb{R}^{n}~,
 \overline{u}_{\mathcal{D}} \in \mathrm{Im}\, (\mathrm{S}_{00}^{-1}\boldsymbol{\mathcal{B}}_{\mathcal{D}}) ,
\end{equation}
in the sense that both problems have the same eigenvalues $\mu$'s, and the corresponding eigenvectors are related via
$ \overline{u}_{\mathcal{D}} = \mathrm{S}_{00}^{-1} \boldsymbol{\mathcal{B}}_{\mathcal{D}}^{\sfrac{1}{2}} \overline{\varphi} \in   \mathrm{Im}\, (\mathrm{S}_{00}^{-1}\boldsymbol{\mathcal{B}}_{\mathcal{D}}) $.
\end{lemma} 
This result is also straightforward and can be obtained multiplying \eqref{E:eigprob-3} by $\mathrm{S}_{00}^{-1} \boldsymbol{\mathcal{B}}_{\mathcal{D}} ^{\sfrac{1}{2}}$.

\vspace{5pt}

To that end, establishing the upper and lower bounds for the eigenvalues of \eqref{E:eigprob-4} and due to equivalence of  \eqref{E:eigprob-4}  with \eqref{E:eigprob-3}, and hence \eqref{E:eigprob-2}, we obtain that eigenvalues of \eqref{E:eigprob-1} are bounded. 
We are interested in nonzero eigenvalues of \eqref{E:eigprob-4} 
for which the following result holds. 

\begin{lemma} \label{L:fact4}   
Let the triangulation $\Omega_h$ for \eqref{E:infinite-pde} be conforming and quasi-uniform. 
Then there exists $\hat{\mu}_\star >0$ independent of the mesh size $h>0$  such that 
\begin{equation} \label{E:main-bounds}
\hat{\mu}_\star \leq \frac{(\boldsymbol{\mathcal{B}}_{\mathcal{D}} \, \overline{u}_{\mathcal{D}} ,  \overline{u}_{\mathcal{D}} )}{(\mathrm{S}_{00} \overline{u}_{\mathcal{D}}, \overline{u}_{\mathcal{D}} )} \leq 1, \quad
\mbox{for all} \quad 0\neq \overline{u}_{\mathcal{D}} ~ 
 \in \mathrm{Im} \,(\mathrm{S}_{00}^{-1}\boldsymbol{\mathcal{B}}_{\mathcal{D}}).
\end{equation}
\end{lemma}

\section{Proofs of statements in Chapter \ref{S:main-results}} \label{S:proof}

\subsection{Harmonic extensions}

Hereafter, we will use the index $\mathcal{D}$ to indicate vectors or functions associated with the domain $\mathcal{D}$ that is the union of all inclusions, and index $0$ to indicate quantities that are associated with the domain outside the inclusions $\Omega\setminus \overline{\mathcal{D}}$.

Now we recall some classical results from the theory of elliptic PDEs. Suppose a function $u^\mathcal{D} \in H^1(\mathcal{D})$,
then consider its harmonic extension  $u^0\in H^1(\Omega \setminus \overline{\mathcal{D}})$ that satisfies
\begin{equation} \label{E:extens-prob}
\left\{
\begin{array}{r l l}
-\bigtriangleup u^0 & =  0, &  \mbox{in } \Omega \setminus \overline{\mathcal{D}} , \\[2pt]
u^0& =  u^\mathcal{D}, &   \mbox{on } \partial \mathcal{D}, \\[2pt]
u^0 & =  0, &  \mbox{on } \Gamma.
\end{array}
\right.
\end{equation}
For such functions the following holds true:
\begin{equation} \label{E:min-fnl}
\displaystyle \int \limits_{\Omega} |\nabla u|^2\;\mathrm{d} x = \min_{v \in H^1_0(\Omega)} \displaystyle \int \limits_{\Omega} |\nabla v|^2\;\mathrm{d} x,
\end{equation}
where
\[
u = \begin{cases} 
u^\mathcal{D}, & \text{ in } \mathcal{D} \\
u^0, & \text{ in } \Omega \setminus \overline{\mathcal{D}}
\end{cases}
\qquad\text{and}\qquad
v = \begin{cases} u^\mathcal{D}, & \text{ in } \mathcal{D} \\
v^0, & \text{ in }  \Omega \setminus \overline{ \mathcal{D}}
\end{cases}
\]
where the function $v^0 \in H^1(\Omega \setminus \overline{\mathcal{D}})$ such that $v^0 |_{\Gamma}=0$, 
and 
\begin{equation} \label{E:exten-norm}
\|u \|_{H^1_0(\Omega)} \leq C \|u^\mathcal{D}\|_{H^1(\mathcal{D})} \quad \mbox{ with the constant}~ C ~\mbox{independent of }~u^\mathcal{D},
\end{equation}
where $\| \cdot \|_{H^1(\Omega)}$ denotes the standard norm of $H^1(\Omega)$:
\begin{equation} \label{E:norm-stand}
 \|v\|^2_{H^1(\Omega)} =  \int_{\Omega} |\nabla v|^2dx + \int_{\Omega} v^2 dx,
\end{equation}
and $\displaystyle \|v \|_{H^1_0(\Omega)}^2 =  \int_{\Omega} |\nabla v|^2dx $.

In view of \eqref{E:min-fnl}, the function $u^0$ of \eqref{E:min-fnl} is the {\it best extension} of $u^\mathcal{D} \in H^1(\mathcal{D})$ among all $H^1(\Omega \setminus \overline{\mathcal{D}})$ functions that vanish on $\Gamma$.
The  algebraic linear system that corresponds to \eqref{E:min-fnl} satisfies the similar property. Namely, if the vector $\overline{u}_0\in \mathbb{R}^{n_{0}}$ is a FEM discretization of the function $u^0 \in H^1_0(\Omega\setminus \overline{\mathcal{D}})$ of \eqref{E:extens-prob}, then for a given $\overline{u}_\mathcal{D} \in \mathbb{R}^{n}$, the best extension $\overline{u}_0  \in \mathbb{R}^{n_{0}}$ would satisfy
\begin{equation} \label{E:min-alg}
\mathrm{A}_{0\mathcal{D}} \,\overline{u}_\mathcal{D}  + \mathrm{A}_{00} \,\overline{u}_0 = 0,
\end{equation}
and
\begin{equation} \label{E:min-algebr}
\left (\bold{A} \begin{bmatrix}\overline{u}_\mathcal{D} \\ \overline{u}_0  \end{bmatrix}, \begin{bmatrix} \overline{u}_\mathcal{D} \\ \overline{u}_0  \end{bmatrix} \right) = \min_{ \overline{v}_0 \in \mathbb{R}^{n_{0}}} \left (\bold{A} \begin{bmatrix} \overline{u}_\mathcal{D} \\ \overline{v}_0 \end{bmatrix},\begin{bmatrix} \overline{u}_\mathcal{D} \\ \overline{v}_0 \end{bmatrix} \right ).
\end{equation}

\subsection{Proof of Lemma \ref{L:fact2} }

Consider generalized eigenvalue problem \eqref{E:eigprob-2} and replace $\boldsymbol{\mathcal{B}}_{\mathcal{D}} $ with $\boldsymbol{\mathcal{B}}_{\mathcal{D}} ^{\sfrac{1}{2}}\boldsymbol{\mathcal{B}}_{\mathcal{D}} ^{\sfrac{1}{2}}$ there, then
\[ 
\boldsymbol{\mathcal{B}}_{\mathcal{D}} ^{\sfrac{1}{2}}\boldsymbol{\mathcal{B}}_{\mathcal{D}} ^{\sfrac{1}{2}} \mathrm{S}_{00}^{-1} \boldsymbol{\mathcal{B}}_{\mathcal{D}} ^{\sfrac{1}{2}}\boldsymbol{\mathcal{B}}_{\mathcal{D}} ^{\sfrac{1}{2}} \,\overline{\psi} = \mu \boldsymbol{\mathcal{B}}_{\mathcal{D}} ^{\sfrac{1}{2}}\boldsymbol{\mathcal{B}}_{\mathcal{D}} ^{\sfrac{1}{2}} \,\overline{\psi}.
\]
Now multiply both sides by the Moore-Penrose 
pseudo inverse $\left[\boldsymbol{\mathcal{B}}_{\mathcal{D}} ^{\sfrac{1}{2}} \right]^{\dagger}$:  
\[
\left[\boldsymbol{\mathcal{B}}_{\mathcal{D}}^{\sfrac{1}{2}} \right]^{\dagger} \boldsymbol{\mathcal{B}}_{\mathcal{D}} ^{\sfrac{1}{2}}\boldsymbol{\mathcal{B}}_{\mathcal{D}}^{\sfrac{1}{2}} \mathrm{S}_{00}^{-1} \boldsymbol{\mathcal{B}}_{\mathcal{D}} ^{\sfrac{1}{2}}\boldsymbol{\mathcal{B}}_{\mathcal{D}} ^{\sfrac{1}{2}}\, \overline{\psi} = \mu \left[\boldsymbol{\mathcal{B}}_{\mathcal{D}} ^{\sfrac{1}{2}} \right]^{\dagger} \boldsymbol{\mathcal{B}}_{\mathcal{D}} ^{\sfrac{1}{2}}\boldsymbol{\mathcal{B}}_{\mathcal{D}} ^{\sfrac{1}{2}} \,\overline{\psi}.
\]
This pseudo inverse has the property that 
\[
\left[\boldsymbol{\mathcal{B}}_{\mathcal{D}}^{\sfrac{1}{2}}\right]^{\dagger} \boldsymbol{\mathcal{B}}_{\mathcal{D}} ^{\sfrac{1}{2}} = \mathrm{P}_{\text{im}},
\] 
where $\mathrm{P}_{\text{im}}$ is an orthogonal projector onto the image $\boldsymbol{\mathcal{B}}_{\mathcal{D}}^{\sfrac{1}{2}}$, hence, $ \mathrm{P}_{\text{im}}\boldsymbol{\mathcal{B}}_{\mathcal{D}}^{\sfrac{1}{2}} = \boldsymbol{\mathcal{B}}_{\mathcal{D}}^{\sfrac{1}{2}}$ and therefore,
\[
\boldsymbol{\mathcal{B}}_{\mathcal{D}}^{\sfrac{1}{2}} \mathrm{S}_{00}^{-1} \boldsymbol{\mathcal{B}}_{\mathcal{D}}^{\sfrac{1}{2}} \overline{\varphi} = \mu\overline{\varphi}, 
\quad \mbox{where}\quad \overline{\varphi} = \boldsymbol{\mathcal{B}}_{\mathcal{D}}^{\sfrac{1}{2}}\overline{\psi}.
\]
Conversely, consider the eigenvalue problem \eqref{E:eigprob-3}, and
multiply its both sides by $\boldsymbol{\mathcal{B}}_{\mathcal{D}} ^{\sfrac{1}{2}} $. Then
\[
\boldsymbol{\mathcal{B}}_{\mathcal{D}} ^{\sfrac{1}{2}}  \boldsymbol{\mathcal{B}}_{\mathcal{D}} ^{\sfrac{1}{2}} \mathrm{S}_{00}^{-1} \boldsymbol{\mathcal{B}}_{\mathcal{D}} ^{\sfrac{1}{2}}  \,\overline{\varphi} 
= \mu \boldsymbol{\mathcal{B}}_{\mathcal{D}} ^{\sfrac{1}{2}} \,\overline{\varphi},
\]
where we replace $\overline{\varphi}$ by $\boldsymbol{\mathcal{B}}_{\mathcal{D}} ^{\sfrac{1}{2}} \, \overline{\psi}$
\[ 
\boldsymbol{\mathcal{B}}_{\mathcal{D}} ^{\sfrac{1}{2}}\boldsymbol{\mathcal{B}}_{\mathcal{D}} ^{\sfrac{1}{2}} \mathrm{S}_{00}^{-1} \boldsymbol{\mathcal{B}}_{\mathcal{D}} ^{\sfrac{1}{2}}\boldsymbol{\mathcal{B}}_{\mathcal{D}} ^{\sfrac{1}{2}} \,\overline{\psi} = \mu \boldsymbol{\mathcal{B}}_{\mathcal{D}}^{\sfrac{1}{2}}\boldsymbol{\mathcal{B}}_{\mathcal{D}} ^{\sfrac{1}{2}} \,\overline{\psi} 
\]
to obtain \eqref{E:eigprob-2}. \hspace{5pt}
$\Box$ 

\subsection{Proof of Lemma \ref{L:fact4} }

\noindent {\bf I. Upper Bound for the Generalized Eigenvalues of \eqref{E:eigprob-1} }

\noindent Consider $\overline{u}=\begin{bmatrix} \overline{u}_{\mathcal{D}}  \\ \overline{u}_0\end{bmatrix}\in\mathbb{R}^N$ with  $\overline{u}_{\mathcal{D}}  
 \in \mathrm{Im}\, (\mathrm{S}_{00}^{-1}\boldsymbol{\mathcal{B}}_{\mathcal{D}} )$, satisfying \eqref{E:min-alg}, 
then 
\begin{equation} \label{E:lem3}
\left (\mathrm{S}_{00} \overline{u}_{\mathcal{D}} ,  \overline{u}_{\mathcal{D}} \right ) = \left (\mathbf{A} \overline{u}, \overline{u} \right ).
\end{equation}
Using \eqref{E:A-def} and \eqref{E:B-def} we obtain from \eqref{E:lem3}:
\begin{equation} \label{E:upper-bd}
\mu= \frac{(\boldsymbol{\mathcal{B}}_{\mathcal{D}} \overline{u}_{\mathcal{D}}, \overline{u}_{\mathcal{D}} )}{(\mathrm{S}_{00}\,  \overline{u}_{\mathcal{D}} , \overline{u}_{\mathcal{D}} )} =
\frac{(\boldsymbol{\mathcal{B}}_{\mathcal{D}} \overline{u}_{\mathcal{D}} ,  \overline{u}_{\mathcal{D}} )}{(\mathbf{A}  \overline{u},  \overline{u})}=  \frac{\displaystyle \int \limits_{\mathcal{D}_h} |\nabla u^{\mathcal{D}}_{h}|^2\;\mathrm{d} x}{\displaystyle \int \limits_{\Omega_h} |\nabla u_h|^2\;\mathrm{d} x}\leq 1,
\end{equation}
with 
\begin{equation} \label{E:u_h}
u_h = \begin{cases}
u^{\mathcal{D}}_{h}, & \mbox{in }\mathcal{D}_h\\[2pt]
u_{h}^0, &  \mbox{in } \Omega \setminus \overline{\mathcal{D}}_h
 \end{cases}
\end{equation}
where $u_{h}^0$ is the harmonic extension of $u^{\mathcal{D}}_{h}$ into $\Omega_h \setminus \overline{\mathcal{D}}_h$ in the sense \eqref{E:extens-prob}. 
$\Box$

\vspace{5pt}

\noindent {\bf II. Lower Bound for the Generalized Eigenvalues of \eqref{E:eigprob-1} }

\begin{figure}[ht]
\centering
\includegraphics[scale=.75]{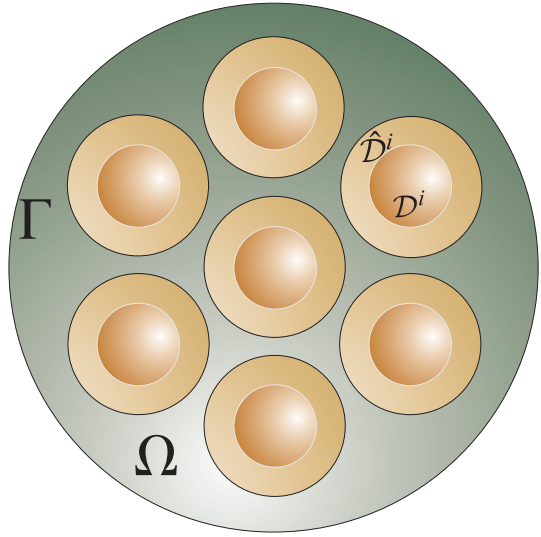} \caption{New domains $\hat{\mathcal{D}}^i$ for our construction of the lower bound of $\mu$} \label{F:new}
\end{figure}

\noindent Before providing the proofs, we introduce one more construction to simplify our consideration below.
Because all inclusions are located at distances that are comparable to their sizes, we construct new domains $\hat{\mathcal{D}}^i$, $i \in \{1,\ldots,m\}$, see Fig. \ref{F:new}, centered at the centers of the original inclusions $\mathcal{D}^i$, $i\in \{1,\ldots,m\}$, but of sizes much larger of those of $\mathcal{D}^i$ and such that
\[
\hat{\mathcal{D}}^i \cap \hat{\mathcal{D}}^j = \emptyset, \quad \mbox{for} \quad i\neq j.
\] 
From it follows below, one can see that the problem \eqref{E:infinite-pde} might be partitioned into $m$ independent subproblems, with what, without loss of generality, we assume that there is only one inclusion, that is, $m=1$.

We also recall a few important results from classical PDE theory analogs of which will be used below. Namely, for a given $v\in H^1(\mathcal{D})$ there exists an extension $v_0$ of $v$ to $\Omega\setminus\overline{\mathcal{D}}$ so that
\begin{equation} \label{E:est-clas}
\| v_0 \|_{H^1(\Omega\setminus\mathcal{D})} \leq C \|v\|_{H^1(\mathcal{D})}, \quad \mbox{with} \quad  C=C(d,\mathcal{D},\Omega).
\end{equation}
One can also introduce a number of norms equivalent to \eqref{E:norm-stand}, and, in particular, below we will use 
\begin{equation} \label{E:norm-our}
 \|v\|^2_{\mathcal{D}}: =  \int_{\mathcal{D}} |\nabla v|^2dx + \frac{1}{R^2}\int_{\mathcal{D}} v^2 dx,
\end{equation}
where $R$ is the radius of the particle $\mathcal{D} = \mathcal{D}_1$. The scaling factor $1/R^2$ is needed for transforming the classical results from a reference (i.e. unit) disk to the disk of radius $R\neq 1$.

We note that the FEM analog of the extension result of \eqref{E:est-clas} for a regular grid was shown in \cite{widlund87}, from which it also follows that the constant $C$ of \eqref{E:est-clas} is independent of the mesh size $h>0$. We utilize this observation in our construction below.

Consider $u_h\in V_h$ given by \eqref{E:u_h}.
Introduce a space $\hat{V}_{h} = \left \{ v_h \in V_{h}:~ v_h = 0  \mbox{ in } \Omega_h \setminus \overline{\hat{\mathcal{D}}}_h \right \}$. Similarly to \eqref{E:u_h}, define
\begin{equation} \label{E:h-hat}
\hat{V}_{h} \ni \hat{u}_h = 
\begin{cases} 
u^{\mathcal{D}}_h,  & \mbox{in }  \mathcal{D}_h\\
\hat{u}^0_h,  &  \mbox{in } \Omega_h  \setminus \overline{\mathcal{D}}_h
\end{cases},
\end{equation}
where $\hat{u}^0_h$ is the harmonic extension of $u^{\mathcal{D}}_h$ into $ \hat{\mathcal{D}}_h \setminus \overline{\mathcal{D}}_h$ in the sense \eqref{E:extens-prob} and $\hat{u}^0_h =0$ on $\partial \hat{\mathcal{D}}_h$.  Also, by \eqref{E:min-fnl} we have
\[
\displaystyle \int \limits_{\Omega_h \setminus \mathcal{D}_h} |\nabla u^0_h|^2 dx \leq \displaystyle \int \limits_{\Omega_h \setminus \mathcal{D}_h} |\nabla \hat{u}^0_h |^2dx.
\]
Define the matrix 
\[
\hat{\mathbf{A}}: = \begin{bmatrix}
\mathrm{A}_{\mathcal{DD}} & \hat{\mathrm{A}}_{\mathcal{D}0} \\
\hat{\mathrm{A}}_{0\mathcal{D}} & \hat{\mathrm{A}}_{00}
\end{bmatrix}
\]
by 
\[
\left ( \hat{\mathbf{A}} \overline{v}, \overline{w}\right ) = \displaystyle \int \limits_{\Omega_h } \nabla {v}_{h} \cdot \nabla {w}_{h} dx, \quad 
\mbox{where} ~ \overline{v}, \overline{w} \in \mathbb{R}^{N}, ~  ~v_h,w_h\in \hat{V}_{h}.
\]
As before, introduce the Schur complement to the block $\hat{\mathrm{A}}_{00}$ of $\hat{\mathbf{A}}$:
\begin{equation} \label{E:Schur-hat}
\hat{\mathrm{S}}_{00} = \mathrm{A}_{\mathcal{DD}} - \hat{\mathrm{A}}_{\mathcal{D}0}\hat{\mathrm{A}}_{00}^{-1}\hat{\mathrm{A}}_{0\mathcal{D}},
\end{equation}
and consider a new generalized eigenvalue problem
\begin{equation} \label{E:eigprob-5}
\boldsymbol{\mathcal{B}}_{\mathcal{D}} \, \overline{u}_{\mathcal{D}}  =  \hat{\mu} \hat{\mathrm{S}}_{00}  \, \overline{u}_{\mathcal{D}} \quad \mbox{with} \quad  \overline{u}_{\mathcal{D}} 
 \in \mathrm{Im}\, (\mathrm{S}_{00}^{-1}\boldsymbol{\mathcal{B}}_{\mathcal{D}}) .
\end{equation}
By \eqref{E:min-algebr} and   \eqref{E:lem3} we have
\begin{equation} \label{E:S-prop}
\left ( \mathrm{S}_{00} \overline{u}_{\mathcal{D}} , \overline{u}_{\mathcal{D}} \right ) \leq \left (\hat{\mathrm{S}}_{00}  \overline{u}_{\mathcal{D}} ,  \overline{u}_{\mathcal{D}} \right) \quad   \mbox{for all }\quad \overline{u}_{\mathcal{D}} 
 \in \mathrm{Im} \,(\mathrm{S}_{00}^{-1}\boldsymbol{\mathcal{B}}_{\mathcal{D}} ).
\end{equation}
Now, we consider a new generalized eigenvalue problem similar to one in \eqref{E:eigprob-3}, namely,
\begin{equation} \label{E:eigprob-6}
\boldsymbol{\mathcal{B}}_{\mathcal{D}}^{\sfrac{1}{2}}  \hat{\mathrm{S}}_{00}^{-1} \boldsymbol{\mathcal{B}}_{\mathcal{D}}^{\sfrac{1}{2}}   \, \overline{\varphi} =   \hat{\mu}\,   \overline{\varphi},  \quad
\overline{\varphi} \in \mathrm{Im} \,  \boldsymbol{\mathcal{B}}_{\mathcal{D}}  .
\end{equation}
We plan to replace $\boldsymbol{\mathcal{B}}_{\mathcal{D}}^{\sfrac{1}{2}}$ in \eqref{E:eigprob-6} with a new symmetric positive-definite matrix $ \hat{\boldsymbol{\mathcal{B}}}_{\mathcal{D}}^{\sfrac{1}{2}}$, given below in \eqref{E:B-hat}, so that 
\begin{equation} \label{E:prop-Bhat}
\boldsymbol{\mathcal{B}}_{\mathcal{D}}^{\sfrac{1}{2}}\boldsymbol{\mathcal{B}}_{\mathcal{D}}^{\sfrac{1}{2}} \overline{\xi}=
\boldsymbol{\mathcal{B}}_{\mathcal{D}}^{\sfrac{1}{2}} \hat{\boldsymbol{\mathcal{B}}}_{\mathcal{D}}^{\sfrac{1}{2}} \overline{\xi}=
 \hat{\boldsymbol{\mathcal{B}}}_{\mathcal{D}}^{\sfrac{1}{2}} \boldsymbol{\mathcal{B}}_{\mathcal{D}}^{\sfrac{1}{2}} \overline{\xi} \quad \mbox{for all} \quad  \overline{\xi}  \in \mathrm{Im} \,  \boldsymbol{\mathcal{B}}_{\mathcal{D}},
\end{equation}
with what \eqref{E:eigprob-6} has the same nonzero eigenvalues as the problem 
\begin{equation} \label{E:eigprob-7}
 \hat{\boldsymbol{\mathcal{B}}}_{\mathcal{D}}^{\sfrac{1}{2}} \hat{\mathrm{S}}_{00}^{-1}  \hat{\boldsymbol{\mathcal{B}}}_{\mathcal{D}}^{\sfrac{1}{2}}  \, \overline{\varphi} =   \hat{\mu} \, \overline{\varphi}, \quad \overline{\varphi}  \in \mathrm{Im} \,  \boldsymbol{\mathcal{B}}_{\mathcal{D}}.
\end{equation}
For this purpose, we consider the decomposition:
\[
\boldsymbol{\mathcal{B}}_{\mathcal{D}} = \mathrm{W} \Lambda \mathrm{W}^T,
\]
where $\mathrm{W}\in \mathbb{R}^{n\times n}$ is an orthogonal matrix composed of eigenvectors $\overline{w}_i$, $i\in \{0,1,\ldots,n-1\}$, of 
\[
\boldsymbol{\mathcal{B}}_{\mathcal{D}} \overline{w} = \nu \overline{w}, \quad \overline{w} \in \mathbb{R}^{n},
\]
and 
\[
\Lambda = \mbox{diag} \left[ \nu_0, \,  \nu_1, \ldots, \nu_{n-1} \right].
\]
Then $\overline{w}_0$ is an eigenvector of $\boldsymbol{\mathcal{B}}_{\mathcal{D}}$ corresponding to $\nu_0=0$ and
\[
\overline{w}_0 = \frac{1}{\sqrt{n}} \begin{bmatrix} 1 \\ \vdots \\ 1 \end{bmatrix}.
\]
To that end, we choose
\begin{equation} \label{E:B-hat}
 \hat{\boldsymbol{\mathcal{B}}}_{\mathcal{D}} = \boldsymbol{\mathcal{B}}_{\mathcal{D}} + \beta \, \overline{w}_0\otimes \overline{w}_0 = \boldsymbol{\mathcal{B}}_{\mathcal{D}} + \beta \, \overline{w}_0 \overline{w}_0^T,
\end{equation}
where $\beta >0 $ is some constant parameter chosen below. Note that the matrix $ \hat{\boldsymbol{\mathcal{B}}}_{\mathcal{D}}$ is symmetric and positive-definite, and  satisfies \eqref{E:prop-Bhat}. It is trivial to show that $ \hat{\boldsymbol{\mathcal{B}}}_{\mathcal{D}}$ given by \eqref{E:B-hat} is spectrally equivalent to $\boldsymbol{\mathcal{B}}_{\mathcal{D}} + \beta \mathrm{I}$ for any $\beta >0$. Also, for quasi-uniform grids, the matrix $h^2 \mathrm{I}$ (in 3-dim case, $h^3 \mathrm{I}$) is spectrally equivalent to the mass matrix $\mathrm{M}_{\mathcal{D}}$ given by
\[
(\mathrm{M}_{\mathcal{D}} \overline{u} , \overline{v} ) = \int_{\mathcal{D}^1_h}  u_h   v_h ~dx, \quad\mbox{where}\quad \overline{u} , \overline{v} \in \mathbb{R}^{n_1}, \quad
u_h , v_h \in V^1_h, 
\]see e.g. \cite{toswid05}. This implies there exists a constant $C>0$ independent of $h$, such that 
\begin{equation} \label{E:spectr-eq}
\left(  \hat{\boldsymbol{\mathcal{B}}}_{\mathcal{D}} \overline{u}_{\mathcal{D}}, \overline{u}_{\mathcal{D}}\right) \geq C   \left( \left(\boldsymbol{\mathcal{B}}_{\mathcal{D}} + \frac{1}{R^2} \mathrm{M}_{\mathcal{D}}\right) \overline{u}_{\mathcal{D}}, \overline{u}_{\mathcal{D}} \right), \quad \mbox{with}\quad \beta=\frac{h^2}{R^2}.
\end{equation}
The choice of the matrix $\boldsymbol{\mathcal{B}}_{\mathcal{D}} + \frac{1}{R^2} \mathrm{M}_{\mathcal{D}}$ for the spectral equivalence was motivated by the fact that 
the right hand side of \eqref{E:spectr-eq} describes $\| \cdot \|_{\mathcal{D}_h}$-norm  \eqref{E:norm-our} of the FEM function $u_h^{\mathcal{D}} \in V_h^1$ that corresponds to the vector $\overline{u}_{\mathcal{D}}\in \mathbb{R}^n$.

Now consider $\overline{u}=\begin{bmatrix} \overline{u}_{\mathcal{D}}  \\ \overline{u}_0\end{bmatrix}\in\mathbb{R}^N$ with  $ \overline{u}_{\mathcal{D}} \in \mathbb{R}^n$, $ \overline{u}_{\mathcal{D}} \in \mathrm{Im}\, \,(\mathrm{S}_{00}^{-1}\boldsymbol{\mathcal{B}}_{\mathcal{D}} )$, and $\overline{u}_0\in\mathbb{R}^{n_0}$ satisfying \eqref{E:min-alg}, and similarly choose 
$\overline{\hat{u}}=\begin{bmatrix} \overline{u}_{\mathcal{D}}  \\ \overline{\hat{u}}_0\end{bmatrix}\in\mathbb{R}^N$ with $\overline{\hat{u}}_0\in\mathbb{R}^{n_0}$ satisfying $\displaystyle \mathrm{\hat{A}}_{0\mathcal{D}} \,\overline{u}_\mathcal{D}  + \mathrm{\hat{A}}_{00} \,\overline{\hat{u}}_0 = 0$, which implies 
\begin{equation} \label{E:S-hat-prop}
\left(\mathrm{\hat{S}}_{00}\overline{u}_\mathcal{D}, \overline{u}_\mathcal{D} \right) = \left( \bold{\hat{A}}\overline{\hat{u}},\overline{\hat{u}}\right). 
\end{equation}
Then 
\begin{equation} \label{E:ext-now}
\left( \bold{\hat{A}}\overline{\hat{u}},\overline{\hat{u}}\right) = \int_{\Omega_h} |\nabla \hat{u}_h|^2dx = 
\int_{\hat{\mathcal{D}}_h\setminus\mathcal{D}_h} |\nabla \hat{u}_h^0|^2dx + \int_{\mathcal{D}_h} |\nabla u_h^\mathcal{D}|^2dx  
\leq (C^*+1) \| u_h^\mathcal{D} \|_{\mathcal{D}_h}^2,
\end{equation}
where $\hat{u}_h \in \hat{V}_h$ is the same extension of $u^{\mathcal{D}}_h$ from $\overline{\mathcal{D}}_h$ to $\Omega_h  \setminus \overline{\mathcal{D}}_h$ as defined in \eqref{E:h-hat}. For the inequality of \eqref{E:ext-now}, we applied the FEM analog of the extension result of \eqref{E:est-clas} by \cite{widlund87}, that yields that the constant $C^*$ in \eqref{E:ext-now} is independent of $h$.

With all the above, we have the following chain of  inequalities:
\[
\begin{array}{l l l}
& \displaystyle 
\frac{\left( \boldsymbol{\mathcal{B}}_{\mathcal{D}} \overline{u}_\mathcal{D}, \overline{u}_\mathcal{D}\right)}{\left (\mathrm{S}_{00}\overline{u}_\mathcal{D}, \overline{u}_\mathcal{D} \right )} 
\underset{\eqref{E:prop-Bhat},\eqref{E:B-hat}}{=}\frac{\left( \left( \boldsymbol{\mathcal{B}}_{\mathcal{D}} + \beta \, \overline{w}_0\otimes \overline{w}_0 \right)\overline{u}_\mathcal{D}, \overline{u}_\mathcal{D}\right)}{\left(\mathrm{S}_{00}\overline{u}_\mathcal{D}, \overline{u}_\mathcal{D} \right)}
\underset{\eqref{E:S-prop}}{\geq}  \frac{\left( \left( \boldsymbol{\mathcal{B}}_{\mathcal{D}} + \beta \, \overline{w}_0\otimes \overline{w}_0 \right)\overline{u}_\mathcal{D},\overline{u}_\mathcal{D}\right)}{\left (\hat{\mathrm{S}}_{00}\overline{u}_\mathcal{D}, \overline{u}_\mathcal{D} \right )} \\[14pt]
 \underset{\eqref{E:S-hat-prop},\eqref{E:spectr-eq}}{\geq} & \displaystyle 
C\frac{\left(  \left( \boldsymbol{\mathcal{B}}_{\mathcal{D}} +\frac{1}{R^2} \mathrm{M}_{\mathcal{D}} \right)\overline{u}_\mathcal{D},\overline{u}_\mathcal{D}\right)}{\left( \bold{\hat{A}}\overline{\hat{u}},\overline{\hat{u}}\right)} \underset{\eqref{E:ext-now}}{\geq} \frac{C \| u^{\mathcal{D}}_h \|_{\mathcal{D}_h}^2 }{(C^*+1)\| u^{\mathcal{D}}_h \|_{\mathcal{D}_h}^2}  =\frac{C }{(C^*+1) } =: \hat{\mu}_\star,
\quad \mbox{with}\quad \beta=\frac{h^2}{R^2}
 \end{array}
\]
where $\mu_\star$ is independent of $h>0$.

\vspace{5pt}

From the obtained above bounds, we have \eqref{E:main-bounds}.
 
$\Box$

\vspace{10pt}

\section{Numerical Results} \label{S:numercs}

In this section, we use four examples to show the numerical
advantages of the Lanczos iterative scheme with the  preconditioner $\boldsymbol{\mathcal{P}}$ defined in \eqref{E:our-prec} 
over the existing preconditioned conjugate gradient method. 

Our numerical experiments are performed by implementing the described above Lanczos algorithm for the problem \eqref{E:pde-form}-\eqref{E:sigma}, where the domain $\Omega$ is chosen to be a disk of radius $5$ with $m=37$ identical circular inclusions $\mathcal{D}^i$, $i\in \{1,\ldots,m\}$. Inclusions are equally spaced. The function $f$ of the right hand side of  \eqref{E:pde-form} is chosen to be a constant, $f=50$. 

In the first set of experiments the values of $\varepsilon_i$'s of \eqref{E:sigma} are going to be identical in all inclusions and vary from $10^{-1}$ to $10^{-8}$.
In the second set of experiments we consider four groups of particles with the same values of $\varepsilon$ in each group that vary from $10^{-4}$ to $10^{-7}$.
In the third set of experiments we consider the case when all inclusions have different values of  $\varepsilon_i$'s that vary from $10^{-1}$ to $10^{-9}$.
Finally, in the fourth set of experiments we decrease the distance between neighboring inclusions.

The initial guess $z^0$ is a random vector that was fixed for all experiments. The stopping criteria is the Euclidian norm of the relative residual $\sfrac{(\boldsymbol{\mathcal{A}}_{\varepsilon} \overline{z}^k-\overline{\mathcal{F}})}{\overline{\mathcal{F}}}$ being less than a fixed tolerance constant.

We test our results agains standard \texttt{pcg} function of $\mbox{MATLAB}^{\text{\textregistered}}$ with $\mathcal{P}_\mathrm{A} = \bold{A}$. The same matrix is also used in the implementation of the described above Lanczos algorithm. 
In the following tables \textbf{PCG} stands for preconditioned conjugate gradient method by $\mbox{MATLAB}^{\text{\textregistered}}$ and \textbf{PL} stands for preconditioned Lanczos method of this paper. 

\textit{Experiment 1.} For the first set of experiments we consider particles $\mathcal{D}^i$ of radius $R=0.45$ in the disk $\Omega$. This choice makes distance $d$ between neighboring inclusions approximately equal to the radius $R$ of inclusions. The triangular mesh $\Omega_{h}$ has $N=32,567$ nodes. Tolerance is chosen to be equal to $10^{-4}$. This experiment concerns the described problem with parameter $\varepsilon$ being the same in each inclusion. Table \ref{table:1} 
shows the number of iterations corresponding to the different values of $\varepsilon$.
\begin{table}[h] \label{T:tab1}
\centering
\caption{Number of iterations  in \textit{Experiment 1}, $N=32,567$}
\begin{tabular}{||c c c c c c c c c||}
 \hline \hline
 & &\\[-2ex]
 $\varepsilon$ &  $10^{-1}$ & $10^{-2}$ & $10^{-3}$ & $10^{-4}$ & $10^{-5}$ & $10^{-6}$ & $10^{-7}$ & $10^{-8}$ \\ [0.5ex] 
 \hline\hline 
 & &\\[-2ex]
 \textbf{PCG} & 10 & 20 & 32 & 40 & 56 & 183 & 302 & 776  \\  [1ex]
 \textbf{PL}  & 33 & 37 & 37 & 37 & 37 & 37 & 37 & 37 \\  [1ex]
\hline \hline
\end{tabular}
\label{table:1}
\end{table}

Based on these results, we first observe that our \textbf{PL} method requires less iterations as $\varepsilon$ goes less than $10^{-4}$. We also notice that number of iterations in the Lanczos algorithm does not depend on $\varepsilon$.

\vspace{5pt}

\begin{figure}[ht]
\centering
\includegraphics[scale=.45]{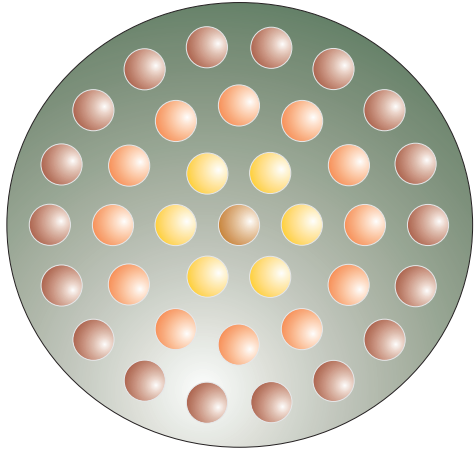} \caption{The domain $\Omega$ with highly conducting inclusions $\mathcal{D}^i$ of fours groups} \label{F:groups}
\end{figure}

\textit{Experiment 2.} In this experiment we leave radii of the inclusions to be the same, namely, $R=0.45$. Tolerance is chosen to be $10^{-6}$. We now distinguish four groups of particles of different $\varepsilon$'s.  The first group consists of one inclusion -- in the center -- with the coefficient $\varepsilon=\varepsilon_1$, whereas the second, third and fourth groups are comprised of the disks in the second, third and fourth circular layers of inclusions with coefficients $\varepsilon_2$, $\varepsilon_3$ and $\varepsilon_4$ respectively, see Fig. \ref{F:groups} (particles of the same group are indicated with the same color). We perform this type of experiments for three different triangular meshes with the total number of nodes $N=5,249$, $N=12,189$ and $N=32,567$. Tables \ref{table:2}, \ref{table:3}, and \ref{table:4} below show the number of iterations corresponding to three meshes respectively.
\begin{table}[h!]
\centering
\caption{Number of iterations  in \textit{Experiment 2}, $N=5,249$}
\begin{tabular}{||c c c c c c||} 
 \hline \hline
 & & & & &\\[-2ex]
 $\varepsilon_1$ & $\varepsilon_2$ &  $\varepsilon_3$ &  $\varepsilon_4$ & \textbf{PCG} & \textbf{PL} \\ [0.5ex] 
 \hline\hline 
 & & & & &\\[-2ex]
$10^{-5}$ & $10^{-5}$ & $10^{-4}$ & $10^{-4}$  & 217 & 39 \\  [1ex]
$10^{-5}$ & $10^{-5}$ & $10^{-4}$ & $10^{-3}$  & 208 & 39 \\  [1ex]
$10^{-6}$ & $10^{-5}$ & $10^{-4}$ & $10^{-3}$  & 716 & 39 \\  [1ex]
$10^{-7}$ & $10^{-6}$ & $10^{-5}$ & $10^{-4}$  & 571 & 39 \\  [1ex]
\hline \hline
\end{tabular}
\label{table:2}
\end{table}
\begin{table}[h!]
\centering
\caption{Number of iterations  in \textit{Experiment 2},  $N=12,189$}
\begin{tabular}{||c c c c c c||} 
 \hline \hline
 & & & & &\\[-2ex]
 $\varepsilon_1$ & $\varepsilon_2$ &  $\varepsilon_3$ &  $\varepsilon_4$ & \textbf{PCG} & \textbf{PL} \\ [0.5ex] 
 \hline\hline 
 & & & & &\\[-2ex]
$10^{-5}$ & $10^{-5}$ & $10^{-4}$ & $10^{-4}$  & 116 & 39 \\  [1ex]
$10^{-5}$ & $10^{-5}$ & $10^{-4}$ & $10^{-3}$  & 208 & 39 \\  [1ex]
$10^{-6}$ & $10^{-5}$ & $10^{-4}$ & $10^{-3}$  & 457 & 39 \\  [1ex]
$10^{-7}$ & $10^{-6}$ & $10^{-5}$ & $10^{-4}$  & 454 & 39 \\  [1ex]
\hline \hline
\end{tabular}
\label{table:3}
\end{table}
\begin{table}[h!]
\centering
\caption{Number of iterations  in \textit{Experiment 2}, $N=32,567$}
\begin{tabular}{||c c c c c c||} 
 \hline \hline
 & & & & &\\[-2ex]
 $\varepsilon_1$ & $\varepsilon_2$ &  $\varepsilon_3$ &  $\varepsilon_4$ & \textbf{PCG} & \textbf{PL} \\ [0.5ex] 
 \hline\hline 
 & & & & &\\[-2ex]
$10^{-5}$ & $10^{-5}$ & $10^{-4}$ & $10^{-4}$  & 311 & 35 \\  [1ex]
$10^{-5}$ & $10^{-5}$ & $10^{-4}$ & $10^{-3}$  & 311 & 35 \\  [1ex]
$10^{-6}$ & $10^{-5}$ & $10^{-4}$ & $10^{-3}$  & 697 & 35 \\  [1ex]
$10^{-7}$ & $10^{-6}$ & $10^{-5}$ & $10^{-4}$  & 693 & 35 \\  [1ex]
\hline \hline
\end{tabular}
\label{table:4}
\end{table}

These results yield that \textbf{PL} requires much less iterations than the corresponding \textbf{PCG} with the number of iterations still being independent of both the contrast $\varepsilon$ and  the mesh size $h$ for \textbf{PL}.

\vspace{5pt}

\textit{Experiment 3.} The next point of interest is to assign different value of $\varepsilon$ for each of 37 inclusions. The geometrical setup is the same as in \textit{Experiment 2}. The value of $\varepsilon_i$, $i\in \{1,\ldots,37\}$, is randomly assigned to each particle and is chosen from the range of $\varepsilon$'s reported in Table \ref{table:5} below. 
The tolerance is $10^{-6}$ as above. The triangular mesh $\Omega_h$ has $12,189$ nodes. 
We run {\bf ten} tests  for each range of contrasts and obtain the {\bf same number of iterations} in every case, and that number is being reported in Table  \ref{table:5}.
\begin{table}[h!]
\centering
\caption{Number of iterations  in \textit{Experiment 3}, $N=12,189$}
\begin{tabular}{||c c ||} 
 \hline \hline
 &\\[-2ex]
Range of $\varepsilon$ & \textbf{PL}  \\ [0.5ex] 
 \hline\hline 
 &\\[-2ex]
$10^{-1}$ to $10^{-8}$  & 53  \\  [1ex]
$10^{-1}$ to $10^{-3}$  & 53 \\  [1ex]
$10^{-7}$ to $10^{-9}$ & 39 \\  [1ex]
\hline \hline
\end{tabular}
\label{table:5}
\end{table}
We also observe that as the contrast between conductivities in the background domain $ \Omega\setminus \overline{\mathcal{D}}$ and the one inside particles $\mathcal{D}_i$, $i\in \{1,\ldots,37\}$, becomes larger our preconditioner demonstrates better convergence, as the third row of Table \ref{table:5} reports. This is expected since the preconditioner constructed above was chosen for the case of absolutely conductive particles. These sets of tests are not compared against the  \textbf{PCG} due to the large number of considered contrasts that prevent this test to converge in a reasonable amount of time.

\textit{Experiment 4.} In the next set of experiments we intend to test how well our algorithm performs if the distance between particles decreases. Recall that the assumption made for our procedure to work is that the interparticle distance $d$ is of order of the particles' radius $R$. With that, we take the same setup as in \textit{Experiment 2} and decrease the distance between particles by making radius of each disk larger. We set $R=0.56$ obtaining that the radius of each inclusion is now twice larger than the distance $d$, and also consider $R=0.59$ so that the radius of an inclusion is three times larger than $d$. The triangular mesh $\Omega_{h}$ has $N=6,329$ and $N=6,497$ nodes, respectively. The tolerance is chosen to be $10^{-6}$. Tables \ref{table:6} and \ref{table:7} show the number of iterations in each case.
\begin{table}[h!]
\centering
\caption{Number of iterations  in \textit{Experiment 4}, $R=0.56$, $N=6,329$}
\begin{tabular}{||c c c c c c||} 
 \hline \hline
 & & & & &\\[-2ex]
 $\varepsilon_1$ & $\varepsilon_2$ &  $\varepsilon_3$ &  $\varepsilon_4$ & \textbf{PCG} & \textbf{PL} \\ [0.5ex] 
 \hline\hline 
 & & & & &\\[-2ex]
$10^{-5}$ & $10^{-5}$ & $10^{-4}$ & $10^{-4}$  & 799 & 61 \\  [1ex]
$10^{-7}$ & $10^{-6}$ & $10^{-5}$ & $10^{-4}$  & 859 & 61 \\  [1ex]
\hline \hline
\end{tabular}
\label{table:6}
\end{table}
\begin{table}[h!]
\centering
\caption{Number of iterations  in \textit{Experiment 4}, $R=0.59$, $N=6,497$}
\begin{tabular}{||c c c c c c||} 
 \hline \hline
 & & & & &\\[-2ex]
 $\varepsilon_1$ & $\varepsilon_2$ &  $\varepsilon_3$ &  $\varepsilon_4$ & \textbf{PCG} & \textbf{PL} \\ [0.5ex] 
 \hline\hline 
 & & & & &\\[-2ex]
$10^{-5}$ & $10^{-5}$ & $10^{-4}$ & $10^{-4}$  & 832 & 73 \\  [1ex]
$10^{-7}$ & $10^{-6}$ & $10^{-5}$ & $10^{-4}$  & 890 & 73 \\  [1ex]
\hline \hline
\end{tabular}
\label{table:7}
\end{table}
Here we observe that number of iterations increases for both \textbf{PCG} and \textbf{PL}, while this number still remains independent of $\varepsilon$ for \textbf{PL}.

We then continue to decrease the distance $d$,  and set $R=0.62$ that is approximately four times larger than the distance between two neighboring inclusions $d$. Choose the same tolerance $10^{-6}$ as above, and the triangular mesh $\Omega_{h}$ of $N=6,699$ nodes, and we observed that our \textbf{PL} method does not reach the desired tolerance in $1,128$ iterations, that confirms our expectations. Further research is needed to develop novel techniques for the case of closely spaced particles that the authors intend to pursue in future.  

\section{Conclusions} \label{S:concl}

This paper focuses on a construction of the robust preconditioner \eqref{E:our-prec} for the Lanczos iterative scheme that can be used in order to solve PDEs with {\it high-contrast} coefficients of the type \eqref{E:pde-form}-\eqref{E:sigma}. A typical FEM discretization yields an ill-conditioning matrix when the contrast in $\sigma$ becomes high (i.e. $\varepsilon \ll 1$). We propose an alternative saddle point formulation \eqref{E:lin-sys-full-matrix-2} of the given problem with the symmetric and indefinite matrix and propose a preconditioner for the employed Lanczos method for solving \eqref{E:lin-sys-full-matrix-2}. The main feature of this novel approach is that we precondition the given linear system with a symmetric and positive semidefinite matrix. The key theoretical outcome is the that the condition number of the constructed preconditioned system is of $O(1)$, which makes the proposed methodology more beneficial for high-contrast problems' application than existing iterative substructuring methods \cite{fr91,flltpr01,mt01,m93,toswid05}. Finally, our numerical results based on simple test scenarios confirm theoretical findings of this paper, and demonstrate convergence of the constructed \textbf{PL} scheme to be independent of the contrast $\varepsilon$, mesh size $h$, and also on the number of different contrasts $\varepsilon_i$, $i\in \{1,\ldots,m\}$ in the inclusions.
In the future, we plan to employ the proposed preconditioner to other types of problems to fully exploit its feature of the independence on contrast and mesh size.

\section{Appendices} \label{S:append}

\subsection{Discussions on the system \eqref{E:lin-sys-full} }   

Along with the problem \eqref{E:lin-sys-full-matrix-1}-\eqref{E:lin-sys-full-matrix-2} and its solution $\bold{x}_\varepsilon$ by \eqref{E:lin-sys-full-matrix-2},
we consider an auxiliary linear system 
\begin{equation} \label{E:lin-sys-infinite-matrix}
\boldsymbol{\mathcal{A}}_o \bold{x}_o=
\begin{bmatrix}
\bold{A} & \bold{B}^T \\
\bold{B} & \bold{0}
\end{bmatrix}
\begin{bmatrix}
\overline{u}_o \\
\overline{\lambda}_o
\end{bmatrix} = 
\begin{bmatrix}
\overline{ \mathrm{F}} \\
\overline{0}
\end{bmatrix},
\end{equation}
or
\begin{equation} \label{E:lin-sys-infinite}
\left\{
\begin{array}{l l l}
\bold{A} \overline{u}_o + \bold{B}^T \overline{\lambda}_o & =\overline{ \mathrm{F}}, \\[2pt]
\bold{B} \overline{u}_o  & =\overline{0}.
\end{array}
\right.
\end{equation}
where matrices $\bold{A}$, $\bold{B}$ and the vector $\overline{ \mathrm{F}}$ are the same as above. 
The linear system \eqref{E:lin-sys-infinite-matrix} or, equivalently \eqref{E:lin-sys-infinite}, emerges in a FEM discretization of the diffusion problem posed in the domain $\Omega$ whose
inclusions are {\it infinitely conducting}, that is, when $\varepsilon=0$ in \eqref{E:sigma}. The corresponding PDE formulation for problem \eqref{E:lin-sys-infinite} might be as follows (see e.g. \cite{cefgal})
\begin{equation} \label{E:infinite-pde}
\left\{
\begin{array}{r l l }
\triangle u & = f, & x \in  \Omega\setminus \overline{\mathcal{D}}  \\[2pt]
u & = \mbox{const}, & x \in \partial \mathcal{D}^i, ~i\in \{1,\ldots,m\} \\[2pt]
\displaystyle \int_{\partial \mathcal{D}^i} \nabla u \cdot \bold{n}_i~ ds & = 0, &  i\in \{1,\ldots,m\} \\[2pt]
u & =0, & x \in  \Gamma
\end{array}
\right.
\end{equation}
where $\bold{n}_i$ is the outer unit normal to the surface $\partial \mathcal{D}^i$. If $u \in H^1_0(\Omega\setminus \overline{\mathcal{D}})$ is an electric potential then it attains constant values on the inclusions $\mathcal{D}^i$ and these constants are not known a priori so that they are unknowns of the problem \eqref{E:infinite-pde} together with $u$.

Formulation  \eqref{E:lin-sys-infinite-matrix} or \eqref{E:lin-sys-infinite} also arises in constrained quadratic optimization problem and solving the Stokes equations for an incompressible fluid \cite{esw05}, and solving elliptic problems using methods combining fictitious domain and distributed Lagrange multiplier techniques to force boundary conditions \cite{gk98}. 

Then the following relation between solutions of systems \eqref{E:lin-sys-full} and \eqref{E:lin-sys-infinite} holds true.
\begin{lemma} \label{L:convergence}
Let $\bold{x}_o=\begin{bmatrix}
\overline{u}_o \\
\overline{\lambda}_o
\end{bmatrix} \in \mathbb{R}^{N+n}$ be the solution of the linear system \eqref{E:lin-sys-infinite}, and 
$\bold{x}_\varepsilon=\begin{bmatrix}
\overline{u} \\
\overline{\lambda}
\end{bmatrix} \in \mathbb{R}^{N+n}$ the solution of \eqref{E:lin-sys-full}. Then 
\[
\overline{u} \to \overline{u}_o \quad \mbox{as} \quad \varepsilon:=\max_{i\in \{1,\ldots,m\}} \varepsilon_i  \to 0.
\]
\end{lemma}

This lemma asserts that the  discrete approximation for the problem \eqref{E:pde-form}-\eqref{E:sigma} converges to the discrete approximation of the solution of \eqref{E:infinite-pde} as $\varepsilon \to 0$. We also note that the continuum version of this fact was shown in \cite{cefgal}.
\begin{proof} 
Hereafter, we denote by $C$ a positive constant that is independent of $\varepsilon$.

Subtract the first equation of  \eqref{E:lin-sys-infinite} from  \eqref{E:lin-sys} and multiply by $\overline{u}_\varepsilon-\overline{u}_o$ to obtain
\[
\left( \bold{A}(\overline{u}_\varepsilon-\overline{u}_o,\overline{u}_\varepsilon-\overline{u}_o) \right ) +
\left( \bold{B}^T(\overline{\lambda}_\varepsilon-\overline{\lambda}_o),\overline{u}_\varepsilon-\overline{u}_o \right ) = \overline{0}.
\]
Recall, that the matrix $\bold{A}$ is SPD then 
\[
(\bold{A}\xi,\xi ) \geq \mu_{1}(\bold{A}) \|\xi\|^2, \quad \forall \xi \in \mathbb{R}^N ,
\] where $\mu_{1}(\bold{A}) >0$ is the {\it minimal eigenvalue} of $\bold{A}$, and $\|\cdot\|=(\cdot,\cdot)$.
Making use of the second equation of \eqref{E:lin-sys-full} we have
\[
\mu_{1}(\bold{A}) \| \overline{u}_\varepsilon-\overline{u}_o \|^2 \leq \left( \boldsymbol{\Sigma}_\varepsilon\overline{\lambda}_\varepsilon, \overline{\lambda}_o \right) - \left( \boldsymbol{\Sigma}_\varepsilon \overline{\lambda}_\varepsilon, \overline{\lambda}_\varepsilon \right) \leq 
 \left(  \boldsymbol{\Sigma}_\varepsilon \overline{\lambda}_\varepsilon, \overline{\lambda}_o \right), 
\]
where we used the fact that $ \boldsymbol{\Sigma}_\varepsilon$ is positive semidefinite. Then
\begin{equation} \label{E:lem1-1}
 \| \overline{u}_\varepsilon-\overline{u}_o \|^2 \leq C\max_i\varepsilon_i \times \| \boldsymbol{\mathcal{B}}_{\mathcal{D}}\overline{\lambda}_\varepsilon \| = C  \varepsilon \| \boldsymbol{\mathcal{B}}_{\mathcal{D}}\overline{\lambda}_\varepsilon \|.
\end{equation}
To estimate $\|\boldsymbol{\mathcal{B}}_{\mathcal{D}}\overline{\lambda}_\varepsilon\|$, we note that $ \boldsymbol{\mathcal{B}}_{\mathcal{D}}\overline{\lambda}_\varepsilon = \bold{B}^T \overline{\lambda}_\varepsilon$, hence, 
\begin{equation} \label{E:lem1-2}
\|  \boldsymbol{\mathcal{B}}_{\mathcal{D}}\overline{\lambda}_\varepsilon \| \leq  \| \overline{ \mathrm{F}} -  \bold{A}\overline{u}_\varepsilon \|.
\end{equation}
Collecting estimates \eqref{E:lem1-1}-\eqref{E:lem1-2}, we observe that it is sufficient to show $\|\overline{u}_\varepsilon\| $ is bounded. For that, we multiply the first equation of \eqref{E:lin-sys} by $\overline{u}_\varepsilon$ and obtain
\[
\left( \bold{A}\overline{u}_\varepsilon,\overline{u}_\varepsilon \right ) +
\left( \bold{B}^T \overline{\lambda}_\varepsilon,\overline{u}_\varepsilon \right ) = \left(\overline{ \mathrm{F}}, \overline{u}_\varepsilon \right),
\quad 
\mbox{that yields}\quad
\| \overline{u}_\varepsilon \| \leq C \| \mathrm{F} \|.
\]
Hence, $ \| \overline{u}_\varepsilon-\overline{u}_o \| \to 0$ as $\varepsilon \to 0$.
\end{proof}

It was also previously observed, see e.g. \cite{irt93,kuz09,toswid05}, that the matrix
\eqref{E:theor-precond}
is the best choice for a preconditioner of $\boldsymbol{\mathcal{A}}_o$. 
This is because there are exactly three  
eigenvalues  of $\boldsymbol{\mathcal{A}}_o$ associated with the following generalized eigenvalue problem 
(see, e.g. \cite{irt93,toswid05})
\begin{equation}   \label{E:gen-evalue-pr}
\boldsymbol{\mathcal{A}}_o \begin{bmatrix}
\overline{u}\\
\overline{\lambda} 
\end{bmatrix} = \mu \bold{P} \begin{bmatrix}
\overline{u}\\
\overline{\lambda}
\end{bmatrix} 
, \quad
\overline{u} \in \mathbb{R}^N, \quad 
\overline{\lambda} \in \mathrm{Im}\,\boldsymbol{\mathcal{B}}_{\mathcal{D}} ,
\end{equation} 
and they are: $\mu_1 < 0$, $\mu_2=1$ and $\mu_3>1$, and, hence, a Krylov subspace iteration method applied for a preconditioned system for solving \eqref{E:gen-evalue-pr} with \eqref{E:theor-precond} {\it converges to the exact solution in three iterations}.

\vspace{5pt}

Now we turn back to the problem \eqref{E:lin-sys-full}. Then the following statement about the generalized eigenvalue problem \eqref{E:gen-eigenpr-eps} holds true.
\begin{lemma} \label{L:indep}
There exist constants $c_1\leq c_2 <0<c_3\leq c_4$ independent of the discretization scale $h>0$ or the contrast parameters $\varepsilon_i$, $i\in \{1,\ldots,m\}$, such that the eigenvalues of the generalized eigenvalue problem
\[
\boldsymbol{\mathcal{A}}_\varepsilon \bold{x}  = \nu \bold{P} \bold{x}  , \quad \overline{u} \in \mathbb{R}^N, ~ 
\overline{\lambda} \in \mathrm{Im}\,\boldsymbol{\mathcal{B}}_{\mathcal{D}},
\]
belong to $[c_1,c_2]\cup[c_3,c_4]$.
\end{lemma}
\noindent Remark that the endpoints $c_i$ of the eigenvalues' intervals might depend on eigenvalues of \eqref{E:gen-evalue-pr}.\\
\noindent \textit{Proof.} \hspace{2pt}
Without loss of generality, here we also assume that all $\varepsilon_i$, $i\in \{1,\ldots,m\}$, are the same and equal to $\varepsilon$, that is, $\bold{\Sigma}_\varepsilon=\varepsilon \boldsymbol{\mathcal{B}}_{\mathcal{D}}$.
Write the given eigenvalue problem as: 
\[
\begin{bmatrix} \bold{A} & \bold{B}^T \\ \bold{B} & - \varepsilon \boldsymbol{\mathcal{B}}_{\mathcal{D}}  \end{bmatrix}
\begin{bmatrix} 
\overline{u}\\
\overline{\lambda} 
\end{bmatrix}
=\nu
 \begin{bmatrix}
\bold{A} & \bold{0}  \\
\bold{0} & \bold{B}\bold{A}^{-1}\bold{B}^T
\end{bmatrix} 
\begin{bmatrix} 
\overline{u}\\
\overline{\lambda} 
\end{bmatrix},
\quad
\overline{u} \in \mathbb{R}^N, ~ 
\overline{\lambda} \in \mathrm{Im}\,\boldsymbol{\mathcal{B}}_{\mathcal{D}} ,\]
which leads to the equation for $\nu$, which is as follows
\begin{equation} \label{E:new-eigenproblem}
\nu - \frac{1}{\nu -1} = -\varepsilon  \, \left[
 \frac{  \left( \boldsymbol{\mathcal{B}}_{\mathcal{D}}\overline{\lambda},\overline{\lambda} \right)  }{  \left(  \bold{B}\bold{A}^{-1}\bold{B}^T \overline{\lambda},\overline{\lambda} \right)  } \right], \quad \overline{\lambda} \in \mathrm{Im}\,\boldsymbol{\mathcal{B}}_{\mathcal{D}}.
\end{equation}
The fraction of the right-hand side of the above equation, that we denote by $\mu$, has been estimated in Theorem \ref{T:main}: $\mu_\star \leq \mu \leq 1$, where $\mu_\star$  is independent of the discretization size $h>0$ due to the norm-preserving extension theorem, \cite{widlund87}.\\
From \eqref{E:new-eigenproblem}, we obtain that the eigenvalues $\nu$ of  \eqref{E:gen-eigenpr-eps} that differ from one, $\nu\neq 1$, are 
\[
\nu^\pm = \frac{1-\varepsilon\mu \pm \sqrt{5 +2 \varepsilon\mu + \varepsilon^2\mu^2 } }{2},
\]
and as $\varepsilon \to 0$, we have $0>\nu^- \to \frac{1- \sqrt{5} }{2}$ and $0<\nu^+ \to \frac{1+ \sqrt{5} }{2}$.
 \\
Finally, using  the bounds for $\mu$ by \eqref{E:main-thm}, we have from \eqref{E:new-eigenproblem} that the endpoints of the intervals $[c_1, c_2]\ni \nu^-$ and $[c_3,c_4] \ni \nu^+$ are independent of both $h$ and $\varepsilon$. In particular, for $0<\varepsilon\ll 1$ and $\mu_\star \leq \mu\leq  1$, we have that $\nu^- \in \left[-\sqrt{2}, \frac{1-\sqrt{5}}{2}\right]$ and $\nu^+ \in \left[\frac{\sqrt{5}}{2}, 2\right]$.  \\[5pt]
If we one has variable $\varepsilon_i$ then it yields a sum over $i\in \{1,\ldots,m\}$ in the right hand side of \eqref{E:new-eigenproblem}. This can be estimated by taking maximal and minimal values of $\varepsilon_i$.\\
$\Box$\\
This lemma demonstrates that \eqref{E:theor-precond} is the best (theoretical) preconditioner for $\boldsymbol{\mathcal{A}}_\varepsilon$ as well as for $\boldsymbol{\mathcal{A}}_o$.

\end{document}